
\baselineskip=14pt
\parskip=10pt

\font\eightrm=cmr8  
\font\eighttt=cmtt8
\magnification=\magstephalf

\parindent=0pt
\overfullrule=0in
\bf
\centerline
{A Comparison Of Two Methods For Random Labelling of 
Balls by Vectors of Integers }
\rm
\bigskip
\centerline{ {\it Doron ZEILBERGER}\footnote{$^1$}
{\eightrm  \raggedright
Department of Mathematics, Temple University,
Philadelphia, PA 19122, USA. 
{\eighttt zeilberg@math.temple.edu .}
Supported in part by the NSF. Version of Dec 19,  1994.
} 
}
 
Greg Kirk[Ki] raised the question of comparing the following two ways
for labelling balls.
Given $r$ pre-determined positive integers $n_i$, $(1 \leq i \leq r)$, and
given $N$  balls ($N$  large), consider two ways  to  randomly assign
$r-$  component  vectors of integers $(a_1  , \dots  ,  a_r)$
to  them,  such  that $1  \leq a_i \leq n_i$. 
We  will call these  vectors  `labels'. Of course  altogether
there are $\prod_{i=1}^r  n_i$ possible labels.
 
{\bf First  Way:}  You put all the balls in one big pot.
For $i=1, \dots  , r$,  at the $i^{th}$ iteration,
line up $n_i$ smaller pots, each with capacity $N/n_i$ balls,
and labeled with labels $1$ through $n_i$,
and, uniformly at random, distribute them into these smaller pots.
Assign the $i^{th}$ component of the vector-label of each ball,
$a_i$, to be the label of the pot in which it was dropped.
Having done that, you dump all the balls back into the big
pot, and go on to the next iteration.
 
{\bf  Second  Way:} Do the same as  above for $i=1$, except that
at the end of the first iteration you do {\it not} dump back the balls
into the large ball but proceed as follows.
For $i=2, \dots , r$, assuming that the balls  have
already   received  their  first  $i-1$  components, leaving the
balls in their pots from the $(i-1)^{th}$ iteration, you line-up
$n_i$ new pots, each with a capacity of $N/n_i$ balls,
and labeled with labels $1$ through $n_i$.
For each of the $n_{i-1}$ pots from the previous iteration, individually,
we uniformly at random, distribute their contents into the new pots,
each of the $n_i$ new pots getting exactly $N/(n_{i-1} n_i)$ balls from
each of the $n_{i-1}$ pots from the previous, $(i-1)^{th}$ iteration.
 
Note that in the First Way, assuming that we can reuse the pots,
we need $1+ max( n_1 ,\dots , n_r)$ pots, one of which should 
have a capacity of $N$ balls, while in the Second Way, we need
$max( 1+ n_1 , n_1+n_2 , \dots , n_{r-1}+n_r )$ pots.
 
The goal is  to maximize the  `equal representation' of all   the possible
$\prod_{i=1}^r n_i$  vector-labels.  It is  obvious, with either
way, that  the probability  of a  ball  to be assigned any  given
label is $\prod_{i=1}^r n_i^{-1}$,  and hence that the  expected number
of balls  to  be given  label $v$, for each of 
$v  \in  \prod_{i=1}^r [1, n_i]$, is $N \prod_{i=1}^r n_i^{-1}$.
 
It is intuitively obvious  that in the  Second Way  the `spread'
in the distribution  is less than in  the First Way. In fact, when
$r=2$, the Second  Way  gives  a  {\it perfect} way of
equi-distribution.  We are {\it guaranteed} that the  number
of  balls given  any particular label $(a_1,a_2)$ is {\it exactly}
$N/  (n_1 n_2)$.
 
Throughout this note  we assume   that $N$ is  divisible by
$lcm( n_1 n_2 , n_2 n_3 , \dots , n_{r-1} n_r )$.
For any  statement $P$, $\chi (P)$ is  $1$ or $0$ according
to whether $P$ is true or false, respectively.
 
The  way to quantify `spread'  is via {\it standard deviation},  or
its square,  the {\it variance}. By symmetry  it is enough
to pick any one {\it fixed} label, $v$, say  $v=(1,  \dots  , 1)$.
 
The `random  variable'  on  a given `experiment'  is the
`number  of balls labeled $v$'. To compute
its variance, we will use an  old  trick, described beautifully
in section 8.2 of the modern classic [GKP]. 
This trick can also be used to  find
the average (i.e.  first  moment), in  which   case  it  is  even  easier
to  use,  and higher moments,  in   which  case it is (usually) harder
to  use.
 
Let  $S$ denote the set of all possible outcomes of the
`labelling  experiment'.
The total number of outcomes, in the
First  Way is

$$
\vert S \vert = \prod_{i=1}^r {{N!} \over  {(N/n_i)!^{n_i}}} \quad .
$$

For each outcome $s$, let $\alpha(s)$   be the quantity
`number of balls that  receive the (fixed) label $v$'.
 
Let's first compute the average of this quantity (even though we
know the answer, just  as a  warm-up for the  calculation  of
the variance, that would follow). We have
 
$$
\sum_{s \in S} \alpha (s)=
\sum_{ s \in S} \sum_{j=1}^{N}  \chi (the \quad j^{th} \quad ball \quad is  
\quad labelled \quad v)
$$
$$
=\sum_{j=1}^N  \sum_{s \in  S_j} 1,   \quad
\eqno(Greg)
$$
 
where the inner sum extends over the set of outcomes, let's call
it $S_j$, of $s \in S$ for which the
$j^{th}$ ball was labelled $v$.
By symmetry, this inner sum is independent  of $j$, and  equals

$$
\prod_{i=1}^r {{(N-1)!} \over  {((N/n_i) -1)! (N/n_i)!^{n_i-1}}} \quad ,
$$
 
since  at each iteration one of the balls (the $j^{th}$) is committed
to lend in one of the pots (Pot $v_i$ in the $i^{th}$ iteration.)
 
Hence the sum in $(Greg)$ equals:

$$
N \prod_{i=1}^r {{(N-1)!} \over  {((N/n_i) -1)! (N/n_i)!^{n_i-1}}} \quad ,
$$
 
and  hence the average is:

$$
av= N \prod_{i=1}^r 
 {
   {
     {
       {(N-1)!} \over  {((N/n_i) -1)! (N/n_i)!^{n_i-1}}   
     } 
   } \over
{
{{(N)!} \over  { (N/n_i)!^{n_i}}} 
}}=
$$
$$
N \prod_{i=1}^r {{1} \over {n_i}}  \quad , 
$$
 
as expected (sic!). 
 
{\bf Variance and Standard Deviation}
 
Let's recall a few elementary facts about {\it variance}. The
{\it standard deviation} is defined to be the  square  root of
the variance. Suppose  that we have a finite set $S$, 
and there is some numerical attribute ({\it random variable}) $X(s)$
for  every element $s \in S$. Then the variance,  $V(X)$, is  the
`average of the  squares  of the `deviation from  the  average'', i.e.
 
$$
V(X)= {{\sum_{s \in S} (X(s)-av)^2}  \over {\vert  S \vert }} \quad  ,
$$
 
where $\vert S  \vert $ is the number of elements  of $S$.
 
It is easier to  compute the related quantity:
 
$$
W(X)= {{\sum_{s \in S} {{X(s)} \choose {2}} }  \over {\vert  S \vert }} 
\quad .
$$
 
Simple algebra shows that:
 
$$
V(X)=2W(X)+av-av^2 \quad.
$$
 
Now we are ready to  compute $W( \alpha )$.
 
$$
 W( \alpha ) =
{{1} \over {\vert S \vert}}  \sum_{ s  \in S } {{\alpha(s)} \choose {2}} =
{{1} \over {\vert S \vert}} \sum_{s  \in S } \sum_{1 \leq i < j \leq N}
\chi( \hbox{ the $i^{th}$ and the $j^{th}$ balls are both labelled $v$})
$$
$$
=
{{1} \over {\vert S \vert}} \sum_{1 \leq i < j  \leq N}
\hbox {[ Number of outcomes with the $i^{th}$ and $j^{th}$ balls both
labelled
$v$]}
\eqno(Kirk)
$$
 
By symmetry, the  summand is independent of  $(i,j)$ and is easily seen to
be equal  to
 
$$
\prod_{i=1}^r 
{{(N-2)!} \over  {((N/n_i) -2)! (N/n_i)!^{n_i-1}}} 
$$
 
since, at each of the $r$ iterations, two balls are committed to   lend
at a predetermined pot (the $v_i^{th}$ pot at the $i^{th}$ iteration.)
 
Simple  algebra  yields
 
$$
W( \alpha )= {{N}  \choose {2}}  
\prod_{i=1}^r n_i^{-2}  
{{(1- n_i/N)} \over {(1- 1/N)}}   \quad .
$$
 
It   follows that
 
$$
V( \alpha)=av-av^2+ 2W(   \alpha)=
{{N} \over {\prod_{i=1}^r n_i  }}-
{{N^2} \over {\prod_{i=1}^r n_i^2  }}+
N(N-1)\prod_{i=1}^r n_i^{-2}  
{{(1- n_i/N)} \over {(1- 1/N)}}  \quad.
$$

Assuming  that $N$ is large, so that  $1/N$  is small,  and using
the  approximation $1/(1-x)=1+x+O(x^2)$, we get the following proposition.
 
{\bf  Proposition 1:}  The  average number of  occurrences of any given vector
$v$ as a  label, in the First Way, is  $N/  \prod_{i=1}^r n_i$,  and its
variance is:
 
$$
{{N} \over { \prod_{i=1}^r n_i } }    -
{{N} \over { \prod_{i=1}^r n_i^2 } } (1+    \sum_{i=1}^r (n_i-1) ) +
O(1) \quad .
$$

{\bf Analysis of the Second Way}
 
Here the total number of outcomes is
 
$$
\vert S \vert =
{{N!} \over {(N/n_1)!^{n_1}}} 
\prod_{i=2}^r 
\left [ {{(N/(n_{i-1})!} \over {(N/(n_{i-1} n_i ) )!^{n_i}}} 
\right ]^{n_{i-1}} \quad.
$$
 
Using an analogous argument as before, the number of outcomes with
the $i^{th}$ and $j^{th}$ balls labelled $v$ equals

$$
{{(N-2)!} \over {((N/n_1) -2)! (N/n_1)!^{n_1-1}}} 
\prod_{i=2}^r 
{{((N/ n_{i-1}) -2)!} \over
{((N/n_{i-1} n_i}) -2)! (N/(n_{i-1} n_i))!^{n_i-1}}
\left [ {{(N/n_{i-1})!} \over {(N/(n_{i-1} n_i ) )!^{n_i} }} \right ]
^{n_{i-1}-1} \quad .
$$
 
Simple algebra yields  that

$$
W( \alpha)=
{{N} \choose {2}} \prod_{i=1}^r  n_i^{-2} \cdot
{{(1- n_1/N)} \over {(1-1/N)}}  \cdot
\prod_{i=2}^r
{{(1- n_{i-1}  n_i/N)} \over {(1- n_{i-1}/N)}}    \quad ,
$$
 
which, as before leads to the following proposition.
 
{\bf  Proposition 2:}  The  average number of  occurrences of any given vector
$v$, as a  label, in the Second Way, is  $N/  \prod_{i=1}^r n_i$,  and its
variance is:
 
$$
{{N} \over { \prod_{i=1}^r n_i } }    -
{{N} \over { \prod_{i=1}^r n_i^2 } } (n_1+    \sum_{i=2}^r (n_i-1)  n_{i-1} ) +
O(1) \quad ,
$$
 
which is slightly smaller.
 
{\bf References}
 
[GKP] R. L. Graham, D. E. Knuth, and O. Patashnik, ``Concrete Mathematics'',
second edition, Addison Wesley, 1993.
 
[Ki] G. Kirk, {\it private communication}, (in `Jeff's Bagels Coffee Bar',
Rocky Hill, NJ, Nov. 1994).
 
\bye